\documentclass[11pt,reqno]{amsart}
\usepackage[a4paper,margin=28mm]{geometry}
\usepackage[T1]{fontenc}
\usepackage{lmodern}
\usepackage{amsmath,amssymb,amsthm,mathtools}
\usepackage{microtype}
\usepackage[ruled,vlined]{algorithm2e}
\usepackage{xcolor}
\usepackage[colorlinks=true,linkcolor=blue!50!black,citecolor=blue!50!black,urlcolor=blue!50!black]{hyperref}
\hypersetup{pdftitle={Dense induced bipartite subgraphs},pdfsubject={Local sparsity, fractional coloring, and induced bipartite subgraphs},pdfauthor={}}
\newtheorem{theorem}{Theorem}[section]
\newtheorem{lemma}[theorem]{Lemma}

\newtheorem{proposition}[theorem]{Proposition}

\newtheorem{conjecture}[theorem]{Conjecture}
\newtheorem{problem}[theorem]{Problem}
\newtheorem{question}[theorem]{Question}
\theoremstyle{remark}
\newtheorem{remark}[theorem]{Remark}
\newcommand{\Pp}{\mathbb P}
\newcommand{\Qq}{\mathbb Q}
\newcommand{\E}{\mathbb E}
\newcommand{\F}{\mathcal F}
\newcommand{\1}{\mathbf 1}
\newcommand{\NL}{N_L}
\newcommand{\NR}{N_R}
\newcommand{\chif}{\chi_f}
\DeclareMathOperator{\Cov}{Cov}

\numberwithin{equation}{section}
\title{Dense induced bipartite subgraphs under local sparsity}
\author[Peiru Kuang, Jing Lin, and Yan Wang]{Peiru Kuang\textsuperscript{1}, Jing Lin\textsuperscript{2}, and Yan Wang\textsuperscript{3}}

\date{}
\begin{document}
\begin{abstract}
The study of dense induced bipartite subgraphs dates back to the work of Erd\H{o}s, Faudree, Pach, and Spencer on triangle-free graphs. In this paper, we study this problem under different local sparsity conditions. 
\begin{itemize}
    \item[(i)] First, we prove that for every fixed integer $t\ge 2$ and all sufficiently large $d$, every $K_{t,t}$-free graph of minimum degree at least $d$ contains an induced bipartite subgraph of average degree at least $\log d-16t\log\log d$. This follows from our bound $\chi_f(G)\le d/(\log d-16t\log\log d)$ for every $d$-degenerate $K_{t,t}$-free graph $G$, which answers a problem of Allen, Dhawan, and Noel.
    \item[(ii)] For triangle-free $d$-degenerate graphs $G$, we prove that $\chi_f(G)\le(2.85+o(1))d/\log d$, improving the leading constant 4 in the bound of Martinsson and Steiner. So every triangle-free graph of minimum degree at least $d$ contains an induced bipartite subgraph of average degree at least $(20/57-o(1))\log d$.
    \item[(iii)] We also extend the study of a conjecture of Esperet, Kang, and Thomass\'e to graphs containing triangles. We prove that if $G$ has minimum degree at least $d$ and each edge belongs to at most $d/k^2$ triangles, where $k=\lfloor\log d/\log\log d\rfloor$, then $G$ contains an induced bipartite subgraph of average degree $\Omega(\log d/\log\log d)$.
    %\item[(iv)] Finally, we give a negative answer to a problem of Allen, Dhawan, and Noel. For infinitely many $d$, we construct a $d$-degenerate graph $G$ with a degeneracy ordering satisfying $|\NL(u)\cap\NR(v)|\le 2\sqrt d$ for all $u,v$, but $\chif(G)\ge(2-o(1))d/\log d$.
\end{itemize}
\end{abstract}
\maketitle

\footnotetext[1]{School of Mathematical Sciences, Shanghai Jiao Tong University, Shanghai 200240, China. Supported by Shanghai Institute for Mathematics and Interdisciplinary Sciences, SIMIS-ID-26-AMS-003. Email: \url{peiru_k@sjtu.edu.cn}.}

\footnotetext[2]{School of Computing and Data Science, Fujian University of Technology, Fujian 350118, China. Supported by Natural
Science Foundation of Fujian Province under Grant No. 2026J001102. Email: \url{sd_frf@163.com}.}

\footnotetext[3]{School of Mathematical Sciences, Shanghai Jiao Tong University, Shanghai 200240, China. Supported by National Key R\&D Program of China under Grant No. 2022YFA1006400 and National Natural Science Foundation of China under Grant No. 12571376. Email: \url{yan.w@sjtu.edu.cn}.}

\section{Introduction}\label{sec:intro}

The \emph{Max-Cut} problem, a central problem in discrete mathematics and theoretical computer
science, seeks to find the largest bipartite subgraph of a graph. For a graph $G$ with $m$ edges, let $f(G)$ denote the maximum number of edges in a bipartite subgraph of $G$. A folklore observation is that $f(G)\ge m/2$ for every graph $G$ with $m$ edges and
the constant $1/2$ in this bound is asymptotically tight~\cite{Erdos1965}. As such, much of the focus has been on determining the error term $f(G)- m/2$. Edwards~\cite{Edwards1973,Edwards1975} proved that
$f(G)-m/2\ge(\sqrt{8m+1}-1)/8$, with equality for complete graphs of odd order. Although the bound $\Omega(\sqrt{m})$ on the error term is optimal in general, it can be significantly improved for
$H$-free graphs, that is, graphs containing no copy of a fixed graph $H$. Dating back to Erd\H{o}s and Lov\'asz (see~\cite{Erdos1979}) in the 1970s, the study of Max-Cut in $H$-free graphs has received significant attention. For example, Shearer~\cite{Shearer1992} proved that every triangle-free graph $G$ with $m$ edges satisfies $f(G)\ge m/2+\Omega(\sum_{v\in V(G)}\sqrt{d_G(v)})$. Alon~\cite{Alon1996Bipartite} subsequently established the sharp bound $f(G)\ge m/2+\Omega(m^{4/5})$. Alon, Krivelevich, and Sudakov~\cite{AKSMaxCut} obtained further bounds for graphs excluding complete bipartite subgraphs. For $d$-degenerate graphs, Carlson, Kolla, Li, Mani, Sudakov, and Trevisan~\cite{CKLMST2021} obtained Max-Cut bounds under restrictions on the number of triangles, while Glock, Janzer, and Sudakov~\cite{GJS2023} established further results for graphs excluding a fixed clique. See also~\cite{ABKS2003,BJS2024} for related results.

The situation, however, is quite different for induced bipartite subgraphs. Indeed, complete graphs have arbitrarily large degree, yet every induced bipartite subgraph of a complete graph has at most two vertices. 
An early motivation for this question was a conjecture of Erd\H{o}s~\cite{Erdos1976}, asserting that every triangle-free graph on $n$ vertices can be made bipartite by deleting at most $n^2/25$ edges. Motivated by this conjecture, Erd\H{o}s, Faudree, Pach, and Spencer~\cite{EFPS} initiated the study of large induced bipartite subgraphs in triangle-free graphs.

An important tool in the study of dense induced bipartite subgraphs is fractional coloring. For a graph $G$, let $\mathcal I(G)$ be the family of its independent sets. A \emph{fractional coloring} assigns nonnegative weights $\lambda_J$ to the sets $J\in\mathcal I(G)$ such that $\sum_{J\ni v}\lambda_J\ge1$ for every vertex $v$. The \emph{fractional chromatic number} $\chif(G)$ is the minimum value of $\sum_{J\in\mathcal I(G)}\lambda_J$. The fractional chromatic number has several equivalent definitions; see~\cite{SU} on fractional coloring as a reference. In this paper, we shall primarily use the following probabilistic characterization: $\chif(G)$ is the minimum real number $q\ge1$ for which there exists a probability distribution on the independent sets of $G$ such that a random independent set $I$ sampled from this distribution satisfies $\Pp(v\in I)\ge1/q$ for every $v\in V(G)$.

The connection between fractional coloring and dense induced bipartite subgraphs is given by the following result of Esperet, Kang, and Thomass\'e~\cite{EKT}: every graph $H$ of average degree $D$ contains an induced bipartite subgraph of average degree at least $D/\chif(H)$. A graph $G$ is \emph{$d$-degenerate} if its vertices admit an ordering $(v_1,\ldots,v_n)$ in which each vertex has at most $d$ neighbors preceding it. To apply our bounds on the fractional chromatic number of $d$-degenerate graphs, we use the following reduction. Every graph of minimum degree at least $d$ contains an induced subgraph that is $d$-degenerate and still has minimum degree at least $d$, see Lemma~\ref{lem:degenerate-reduction}. 

In this paper, we study dense induced bipartite subgraphs under three different local sparsity conditions. The first concerns $K_{t,t}$-free graphs for every integer $t\ge 2$; the second concerns triangle-free graphs; and the third concerns graphs in which adjacent vertices have few common neighbors.

\subsection{\texorpdfstring{$K_{t,t}$}{K(t,t)}-free graphs}

We begin with the following result for $K_{t,t}$-free graphs.

\begin{theorem}\label{thm:bipartite}
For every integer $t\ge2$ and all sufficiently large integers $d$, every $K_{t,t}$-free graph of minimum degree at least $d$ contains an induced bipartite subgraph of average degree at least $\log d-16t\log\log d$.
\end{theorem}

%An upper bound on the fractional chromatic number also yields a lower bound on the independence number, since $\alpha(G)\ge |V(G)|/\chif(G)$. Closely related developments have taken place in graph coloring, where local sparsity can lead to substantial improvements in coloring bounds. Alon, Krivelevich, and Sudakov~\cite{AKS} proved that if $G$ has maximum degree $\Delta$ and every neighborhood spans at most $\Delta^2/f$ edges, then $\chi(G)=O(\Delta/\log f)$. For triangle-free graphs, Molloy~\cite{Molloy} obtained the asymptotically sharper bound $\chi_\ell(G)\le(1+o(1))\Delta/\log\Delta$ for the list chromatic number. More recently, Anderson, Bernshteyn, and Dhawan~\cite{ABD} proved that every $K_{t,t}$-free graph of maximum degree $\Delta$ satisfies $\chi(G)\le(1+o(1))\Delta/\log\Delta$ for fixed $t\ge2$. Davies, Kang, Pirot, and Sereni~\cite{DKPS} developed the local occupancy method, which provides a unified approach to coloring locally sparse graphs. These results motivate the study of coloring bounds in terms of degeneracy rather than maximum degree.

Local sparsity has long been used to obtain improved coloring bounds in terms of maximum degree; see~\cite{AKS,DKPS,Molloy}. In particular, Anderson, Bernshteyn, and Dhawan~\cite{ABD} proved that every $K_{t,t}$-free graph of maximum degree $\Delta$ satisfies $\chi(G)\le(1+o(1))\Delta/\log\Delta$ for each fixed $t\ge2$. However, a $d$-degenerate graph may have arbitrarily large maximum degree, so these results do not directly yield analogous bounds in terms of degeneracy.

Dhawan~\cite{Dhawan} developed an entropy method for sampling random independent sets. Building on this approach, Allen, Dhawan, and Noel~\cite{ADN} proved that $\chif(G)\le(1+o(1))d/\log d$ for every $d$-degenerate $C_4$-free graph $G$. They posed the following problem.

\begin{problem}[Allen--Dhawan--Noel]\label{ques:ADN}
Let $G$ be a $d$-degenerate $K_{t,t}$-free graph. Is $\chi_f(G)\leq(1+o(1))\frac{d}{\log d}$?
\end{problem}

Dhawan, Nguyen, and Rathi~\cite{DNR} recently proved that if a $d$-degenerate graph admits an ordering in which the earlier neighborhood of every vertex spans at most $d^2/f$ edges, then its fractional chromatic number is at most $(8+o(1))d/\log f$. By the K\H{o}v\'ari--S\'os--Tur\'an theorem~\cite{KST}, their result implies that every $d$-degenerate $K_{t,t}$-free graph satisfies $\chif(G)\le(8t+o(1))d/\log d$. 
Bohman and Keevash~\cite{BK2010} used the random $K_{t,t}$-free process to obtain improved lower bounds on $\operatorname{ex}(n,K_{t,t})$ for $t\ge5$.

We answer Problem \ref{ques:ADN} affirmatively in the following quantitative form.

\begin{theorem}\label{thm:main}
For every integer $t\ge2$, there exists $d_0(t)$ such that, for every integer $d\ge d_0(t)$, every $d$-degenerate $K_{t,t}$-free graph $G$ satisfies
\begin{equation}\label{eq:main}
  \chif(G)\le\frac{d}{\log d-16t\log\log d}.
\end{equation}
\end{theorem}

In particular, $\chif(G)\le(1+o_t(1))d/\log d$. The leading constant is best possible: Allen, Dhawan, and Noel~\cite[Theorem~1.4]{ADN} constructed $d$-degenerate graphs of arbitrarily large fixed girth with fractional chromatic number at least $(1-o(1))d/\log d$. 

Theorem~\ref{thm:bipartite} follows directly from Theorem~\ref{thm:main} and the reduction of~\cite{EKT}, after passing to a $d$-degenerate induced subgraph of minimum degree at least $d$ (Lemma \ref{lem:degenerate-reduction}). 
To prove Theorem~\ref{thm:main}, we use a random independent-set construction of Allen, Dhawan, and Noel~\cite{ADN}, but develop a different analysis that applies to $K_{t,t}$-free graphs. Our analysis uses a change of probability measure weighted by the final weight of a fixed vertex. Under this measure, a weighted K\H{o}v\'ari--S\'os--Tur\'an estimate and martingale arguments control the total weight of its earlier neighbors.

\subsection{Triangle-free graphs}

We first give a lower bound on the average degree of an induced bipartite subgraph in a triangle-free graph of large minimum degree.

\begin{theorem}\label{thm:triangle-free-bipartite}
Every triangle-free graph of minimum degree at least $d$ contains an induced bipartite subgraph of average degree at least $(20/57-o(1))\log d$.
\end{theorem}

As in the $K_{t,t}$-free case, we deduce this result from an upper bound on the fractional chromatic number of $d$-degenerate triangle-free graphs. Harris~\cite[Conjecture~6.2]{Harris} conjectured that the trivial bound $d+1$ can be improved by a logarithmic factor, namely that $\chif(G)=O(d/\log d)$ for every triangle-free $d$-degenerate graph $G$; see also~\cite[Conjecture~1.6]{MS}. Martinsson and Steiner~\cite{MS} confirmed this conjecture with the following bound.

\begin{theorem}[Martinsson--Steiner]\label{thm:MS-degenerate}
Suppose $G$ is a triangle-free and $d$-degenerate graph. Then $\chi_f(G)\le(4+o(1))\frac{d}{\log d}$, where the $o(1)$ term tends to $0$ as $d$ increases.
\end{theorem}

They further proposed the following conjecture on the leading constant.

\begin{conjecture}[Martinsson--Steiner]\label{conj:MS-sharp}
The following holds for any sufficiently large $d$.
\begin{enumerate}
\item[1.] $\chi_f(G)\le(1+o(1))\frac{d}{\log d}$ for all $d$-degenerate triangle-free graphs $G$.
\item[2.] There exists a $d$-degenerate triangle-free graph $G$ such that $\chi_f(G)\ge(1-o(1))\frac{d}{\log d}$.
\end{enumerate}
\end{conjecture}

The second part of Conjecture~\ref{conj:MS-sharp} was proved by Allen, Dhawan, and Noel~\cite[Theorem~1.4]{ADN}. We improve the upper bound in Theorem~\ref{thm:MS-degenerate} by reducing the leading constant from $4$ to $2.85$.

\begin{theorem}\label{thm:triangle-free}
Every triangle-free $d$-degenerate graph $G$ satisfies
\begin{equation}\label{eq:triangle-free}
\chif(G)\le\left(2.85+o(1)\right)\frac d{\log d}
\end{equation}
as $d\to\infty$, where the error term is uniform over all such graphs $G$.
\end{theorem}

Our proof is based on two complementary estimates for the same weighted expression. First, we obtain a sharper entropy estimate by repeatedly applying an inequality arising from the procedure of Martinsson and Steiner~\cite[Section~4]{MS}. Second, we introduce a stopping rule for the sampling procedure in Section~\ref{sec:sampler} and derive a lower bound using a submartingale. Comparing these estimates allows us to reduce the leading constant below $3$.

\subsection{Graphs with a few triangles}

The preceding results concern graphs in which certain subgraphs are forbidden. We now consider whether dense induced bipartite subgraphs can still be guaranteed when triangles are allowed. This question is related to a broader problem motivated by separation choosability, a variant of list coloring introduced by Kratochv\'il, Tuza, and Voigt~\cite{KTV1998}, in which the lists assigned to adjacent vertices have at most one color in common. Esperet, Kang, and Thomass\'e~\cite{EKT} proposed the following conjecture.

\begin{conjecture}[Esperet--Kang--Thomass\'e]\label{conj:EKT}
\leavevmode\\[-0.75\baselineskip]
\begin{enumerate}
\item[(i)] Any graph with minimum degree at least $d$ contains either a complete subgraph on $\omega_{d\to\infty}(1)$ vertices or an induced bipartite subgraph with minimum degree $\omega_{d\to\infty}(1)$.
\item[(ii)] There is a constant $c>0$ such that any triangle-free graph with minimum degree at least $d$ contains an induced bipartite subgraph of minimum degree at least $c\log d$.
\end{enumerate}
\end{conjecture}

Here, $\omega_{d\to\infty}(1)$ denotes a function of $d$ that tends
to infinity as $d\to\infty$.
Note that (ii) strengthens the triangle-free case of (i) by providing a logarithmic lower bound on the minimum degree. Kwan, Letzter, Sudakov, and Tran~\cite{KLST} showed that, for every fixed graph $H$, an $H$-free graph of minimum degree at least $d$ contains an induced bipartite subgraph of minimum degree $\Omega_H(\log d/\log\log d)$. Their result confirmed (i) and established (ii) up to a factor of $\log\log d$. Glock~\cite{Glock} subsequently gave a shorter proof of the latter bound for triangle-free graphs. More recently, Martinsson and Steiner~\cite{MS} established (ii), using random independent sets to obtain the required logarithmic lower bound.

We next investigate whether a similar lower bound can be obtained when triangles are allowed. For an edge $uv\in E(G)$, define its \emph{edge-codegree} by $c_G(uv)=|N_G(u)\cap N_G(v)|$, the number of triangles containing $uv$. Our next result shows that the bound $\Omega(\log d/\log\log d)$ continues to hold when every edge belongs to relatively few triangles. For each integer $d\ge3$, let $k(d)=\lfloor\log d/\log\log d\rfloor$.

\begin{theorem}\label{thm:edge-codegree}
There is an absolute constant $c>0$ such that, for every integer $d\ge3$, every graph $G$ with $\delta(G)\ge d$ and $c_G(uv)\le d/k(d)^2$ for all $uv\in E(G)$ contains an induced bipartite subgraph of average degree at least $c\log d/\log\log d$.
\end{theorem}

Since every graph of average degree at least $D$ contains an induced subgraph of minimum degree at least $D/2$, this theorem also guarantees an induced bipartite subgraph of minimum degree $\Omega(\log d/\log\log d)$.

The proof uses a different sampling method. We choose a random subset of vertices and find two disjoint independent sets with sufficiently many edges between them.

\medskip
\noindent
\textbf{Organization.}
The rest of the paper is organized as follows. Section~\ref{sec:prelim} presents the basic definitions and preliminary lemmas. Sections~\ref{sec:sampler} and~\ref{sec:fractional-proof} develop our construction of random independent sets and establish Theorem~\ref{thm:main}, from which Theorem~\ref{thm:bipartite} follows. Section~\ref{sec:triangle-free} proves Theorem~\ref{thm:triangle-free}. Section~\ref{sec:edge} proves Theorem~\ref{thm:edge-codegree} using a different sampling method. Finally, Section~\ref{sec:questions} discusses further questions.

\section{Preliminaries}\label{sec:prelim}
All graphs considered in this paper are finite and simple. For a graph $G$, let $\overline{d}(G)=2|E(G)|/|V(G)|$ denote its average degree. For a subset $U \subseteq V(G)$, we denote by $G[U]$ the subgraph of $G$ induced by $U$. For a vertex $v \in V(G)$, let $N_G(v)$ denote its neighborhood in $G$, or simply $N(v)$ when $G$ is clear from the context. For a real number $x$, let $x_+ = \max\{x, 0\}$ denote its positive part. Throughout this paper, all logarithms are natural. %For an integer $i \ge 0$, the falling factorial is defined by $(a)_i = a(a-1)\cdots(a-i+1)$, with the convention that $(a)_0 = 1$.

%Let $(\F_i)_{i=0}^m$ be a filtration. A sequence $(U_i)_{i=0}^m$ is a \emph{martingale} under $\Pp$ if each $U_i$ is integrable and $\F_i$-measurable, and $\E_{\Pp}[U_i\mid\F_{i-1}]=U_{i-1}$ almost surely for every $1\le i\le m$. It is a \emph{submartingale} if the equality is replaced by $\E_{\Pp}[U_i\mid\F_{i-1}]\ge U_{i-1}$. A random variable $T$ taking values in $\{0,\ldots,m\}$ is a \emph{stopping time} if $\{T\le i\}\in\F_i$ for every $0\le i\le m$.

We use the standard notions of martingales, submartingales, and stopping times; see Williams~\cite[Chapter~10]{Williams}.
For square-integrable random variables $U$ and $V$, we write
\begin{equation}\label{cov}
\Cov_{\Pp}(U,V\mid\F_i)
=\E_{\Pp}[UV\mid\F_i]
-\E_{\Pp}[U\mid\F_i]\E_{\Pp}[V\mid\F_i].
\end{equation}
%For further background on martingales and stopping times, we refer the reader to Williams~\cite{Williams}.

For a fixed $d$-degeneracy ordering, we write $\NL(v_i)=N(v_i)\cap\{v_1,\ldots,v_{i-1}\}$ and $\NR(v_i)=N(v_i)\cap\{v_{i+1},\ldots,v_n\}$. Thus $|\NL(v_i)|\le d$, and in particular $e(G)\le d|V(G)|$. The following standard observation allows us to pass to a $d$-degenerate induced subgraph while preserving the lower bound of minimum degree.

\begin{lemma}\label{lem:degenerate-reduction}
Let $d$ be a positive integer and let $G$ have minimum degree at least $d$. Then $G$ has an induced subgraph $H$ that is $d$-degenerate and has minimum degree at least $d$.
\end{lemma}

\begin{proof}
Choose an induced subgraph $H$ of $G$ with minimum degree at least $d$ and the fewest vertices. By minimality, $\delta(H)=d$, while every proper nonempty induced subgraph of $H$ has minimum degree less than $d$. Hence $H$ is $d$-degenerate.
\end{proof}

The following lemma of Esperet, Kang, and Thomass\'e~\cite[Theorem~3.1]{EKT} allows us to find dense induced bipartite subgraphs in graphs with small fractional chromatic number. It will be used in the proofs of Theorems~\ref{thm:bipartite} and~\ref{thm:triangle-free-bipartite}.

\begin{lemma}[Esperet--Kang--Thomass\'e]\label{lem:bipartite-transfer}
Every nonempty graph $H$ contains an induced bipartite subgraph of average degree at least $\overline d(H)/\chif(H)$.
\end{lemma}

We next recall the following weighted version of the K\H{o}v\'ari--S\'os--Tur\'an theorem due to Fox, Kwan, and Spink~\cite[Lemma~3.3]{FKS}.

\begin{lemma}[Fox--Kwan--Spink]\label{lem:FKS-weighted}
Consider a bipartite graph $G=(A,B)$ where the vertices in $B$ are weighted (with total weight $1$) and $|A|=n$. Let $w(G)=\sum_{a\sim b}w(b)$. Assume every set of $t$ vertices in $A$ has common neighborhood with total weight at most $q$. Then
\[
w(G)\le q^{1/t}n+t-1.
\]
\end{lemma}

The following lemma is a direct consequence of Lemma \ref{lem:FKS-weighted}.

\begin{lemma}\label{lem:weighted}
Let $s\ge2$, $h>0$, and $\delta>0$. Let $G$ be $K_{s+1,s+1}$-free, and let $A\subseteq N(v)$ for a vertex $v$. Assign weights $r_a\in[0,h]$ to the vertices $a\in A$, and write $M=\sum_{a\in A}r_a$. 
If $\delta\ge 2s^{1/s}h^{1/s}M^{1-1/s}$, then for any $U\subseteq V(G)\setminus\{v\}$, we have
\[
\left|\left\{u\in U:\sum_{a\in A\cap N(u)}r_a>\delta\right\}\right|
\le \frac{2(s-1)M}{\delta}.
\]
\end{lemma}

\begin{proof}
Let $W$ be the set of vertices $u\in U$ satisfying $\sum_{a\in A\cap N(u)}r_a>\delta$, and write $m=|W|$. We may assume that $M>0$, since otherwise $W$ is empty.

We first observe that any $s$ vertices of $W$ have at most $s$ common neighbors in $A$. Indeed, if they had $s+1$ common neighbors in $A$, then these vertices, together with $v$, would form a copy of $K_{s+1,s+1}$, a contradiction to our assumption.

Let $B$ be the bipartite graph whose two vertex classes are disjoint copies of $W$ and $A$, where $u\in W$ and $a\in A$ are adjacent if $ua\in E(G)$. Assign each $a\in A$ the weight $r_a/M$. Since any $s$ vertices of $W$ have at most $s$ common neighbors in $A$, their total weight is at most $sh/M$. Applying Lemma~\ref{lem:FKS-weighted} with $q=sh/M$, we obtain
$
\frac{1}{M}\sum_{u\in W}\sum_{a\in A\cap N(u)}r_a
\le \left(\frac{sh}{M}\right)^{1/s}m+s-1.
$
Since $\sum_{a\in A\cap N(u)}r_a>\delta$ for every $u\in W$, it follows that
$
m\delta\le s^{1/s}h^{1/s}mM^{1-1/s}+(s-1)M.
$
By the assumption on $\delta$, the first term on the right is at most $m\delta/2$. It follows that
$
\frac{m\delta}{2}\le(s-1)M,
$
which gives the desired bound.
\end{proof}

\section{Independent set sampling}\label{sec:sampler}
In this section, we describe the random independent set sampling procedure used to prove Theorem~\ref{thm:main}. We then introduce a change of measure to control the probability that a vertex's weight reaches the upper threshold.

\subsection{The sampling procedure}
We use the random independent set procedure of Allen, Dhawan, and Noel~\cite[Algorithm~1]{ADN}, with parameters $0<\alpha<h<1/2$ to be specified later. While their analysis relies on the $C_4$-free assumption, we develop a different approach that applies to $K_{t,t}$-free graphs. Fix a $d$-degeneracy ordering $v_1,\ldots,v_n$ of $G$. For each vertex $w$, let $p_i(w)$ denote its weight after step $i$. Let $\F_i$ be the $\sigma$-algebra generated by all random choices made during the first $i$ steps, with $\F_0$ trivial. 

\medskip
\begin{algorithm}[H]
\caption{Random independent set procedure}\label{alg:random-independent-set}
\KwIn{A graph $G$, an ordering $v_1,\ldots,v_n$, and parameters $0<\alpha<h<1/2$.}
\KwOut{An independent set $I\subseteq V(G)$.}
$I\gets\varnothing$ and $p_0(w)\gets\alpha$ for every $w\in V(G)$\;
\For{$i=1,\ldots,n$}{
  \eIf{$p_{i-1}(v_i)\in\{0,h\}$}{
    $q_i\gets0$, $A_i\gets0$; leave all weights unchanged\;
  }{
    $q_i\gets p_{i-1}(v_i)$\;
    Choose $A_i\in\{0,1\}$ with $\Pp(A_i=1\mid\F_{i-1})=q_i$\;
    \If{$A_i=1$}{Add $v_i$ to $I$\;}
    \ForEach{$w\in\NR(v_i)$ with $p_{i-1}(w)<h$}{
      $x\gets p_{i-1}(w)$, $c\gets\min\{x/(1-q_i),h\}$\;
      $\theta\gets\dfrac{x-(1-q_i)c}{q_i h}$\;
      \eIf{$A_i=0$}{
        $p_i(w)\gets c$\;
      }{
        Independently set $p_i(w)\gets h$ with probability $\theta$, and $p_i(w)\gets0$ otherwise\;
      }
    }
    Leave all other weights unchanged\;
  }
}
\end{algorithm}

\medskip

Conditional on $\F_{i-1}$ and $A_i=1$, the auxiliary random choices for distinct right-neighbors of $v_i$ are independent. We first verify that $\theta\in[0,1]$. If $x/(1-q_i)\le h$, then $c=x/(1-q_i)$ and hence $\theta=0$. Otherwise, $c=h$ and
$
\theta=\frac{x-(1-q_i)h}{q_i h}
=\frac{x/h-1+q_i}{q_i}.
$
In this case, $\theta>0$, and $\theta\le1$ since $x<h$. Thus, all weights remain in $[0,h]$. Moreover, once a weight reaches $0$ or $h$, it remains unchanged. Whenever $v_i$ is selected, each of its right-neighbors has weight $0$ or $h$ after step $i$ and cannot be selected later. Thus, $I$ is an independent set.

By Algorithm \ref{alg:random-independent-set}, we have the following two properties. First, for every vertex $w$ and every $1\le i\le n$,
\begin{equation}\label{eq:weight-growth}
p_i(w)\le\frac{p_{i-1}(w)}{1-q_i}.
\end{equation}
To see this, let $w\in N_R(v_i)$ with $p_{i-1}(w)<h$, and write $x=p_{i-1}(w)$. If $A_i=0$, then $p_i(w)=c\le x/(1-q_i)$. If $A_i=1$, then $p_i(w)$ is either $0$ or $h$, and the latter is possible only when $c=h\le x/(1-q_i)$. The inequality also holds for weights that remain unchanged.

Second, for every vertex $w$ and every $1\le i\le n$,
\begin{equation}\label{eq:original-martingale}
\E_{\Pp}[p_i(w)\mid\F_{i-1}]=p_{i-1}(w).
\end{equation}
Indeed, for $w\in N_R(v_i)$ with $p_{i-1}(w)<h$, write $x=p_{i-1}(w)$. If $q_i>0$, then the conditional expectation of $p_i(w)$ is $(1-q_i)c+q_i h\theta=x$. All other weights remain unchanged. Thus, for each vertex $w$, the process $(p_i(w))_{i=0}^n$ is a martingale under $\Pp$. In particular, $\E_{\Pp}[p_i(w)]=p_0(w)=\alpha$ for every $0\le i\le n$.

For each vertex $v_k$, let $\mathcal B_k$ denote the event $p_{k-1}(v_k)=h$. By above, $\Pp(v_k\in I\mid\F_{k-1})=p_{k-1}(v_k)-h\mathbf1_{\mathcal B_k}$. Taking expectations and using~\eqref{eq:original-martingale} gives
\begin{equation}\label{eq:marginal}
\Pp(v_k\in I)=\E_{\Pp}[p_{k-1}(v_k)]-h\Pp(\mathcal B_k)=\alpha-h\Pp(\mathcal B_k).
\end{equation}

The main difficulty is to show that the weight of a vertex is unlikely to reach the upper threshold, since a vertex can no longer be selected once this happens. Allen, Dhawan, and Noel~\cite{ADN} address this difficulty by bounding the second moment of the vertex's weight using a supermartingale. We take a different approach by studying the total weight of its left-neighbors. The key observation is that a vertex can reach the upper threshold only if this total weight becomes sufficiently large. We use a change of measure to relate the probability of reaching the threshold to the growth of the left-neighbor weights. Under the new measure, Lemma~\ref{lem:FKS-weighted} bounds the number of steps with large contributions to this growth, while martingale estimates control the remaining steps.

\subsection{A change of measure at one vertex}

Fix a vertex $v=v_k$, and let $A=N_L(v)$ and $K=k-1$. If $A=\varnothing$, then $p_K(v)=\alpha$ and hence $\Pp(v\in I)=\alpha$. We may therefore assume that $A\ne\varnothing$.

By~\eqref{eq:marginal}, to show that $\Pp(v\in I)\ge(1-o(1))\alpha$, we need to prove that $h\Pp(\mathcal B_k)/\alpha=o(1)$. Since $h/\alpha$ can be large, it is not enough to show that $\Pp(\mathcal B_k)=o(1)$. To overcome this difficulty, we introduce a new probability measure $\Qq$ under which $h\Pp(\mathcal B_k)/\alpha$ is exactly the probability of $\mathcal B_k$.

Since $p_K(v)=h$ on $\mathcal B_k$ and $\E_{\Pp}[p_K(v)]=\alpha$, we define a probability measure $\Qq$ on $\F_K$ by
\begin{equation}\label{eq:change-measure}
\frac{d\Qq}{d\Pp}=\frac{p_K(v)}{\alpha}.
\end{equation}
Then $\Qq(\mathcal B_k)=h\Pp(\mathcal B_k)/\alpha$, so~\eqref{eq:marginal} becomes
\begin{equation}\label{eq:selection-Q}
\Pp(v\in I)=\alpha(1-\Qq(\mathcal B_k)).
\end{equation}
Thus, it suffices to show that $\Qq(\mathcal B_k)=o(1)$. In fact, we shall prove the stronger bound $\Qq(\mathcal B_k)=O((\log d)^{-2})$.
By~\eqref{eq:original-martingale}, $\E_{\Pp}[p_K(v)\mid\F_i]=p_i(v)$ for every $0\le i\le K$. Hence
\begin{equation}\label{eq:tilted-density}
\frac{d(\Qq|_{\F_i})}{d(\Pp|_{\F_i})}=\frac{p_i(v)}{\alpha}.
\end{equation}
For every $1\le i\le K$ and every bounded $\F_i$-measurable random variable $Z$, the conditional Bayes formula~\cite[Theorem~3.22]{Xiong} gives
\begin{equation}\label{eq:tilted-transition}
\E_{\Qq}[Z\mid\F_{i-1}]
=\frac{\E_{\Pp}[p_i(v)Z\mid\F_{i-1}]}{p_{i-1}(v)}.
\end{equation}
These equalities hold $\Qq$-almost surely, since~\eqref{eq:tilted-density} gives $\Qq(p_{i-1}(v)=0)=0$.

For each $0\le i\le K$, let $S_i=\sum_{a\in A}p_i(a)$. In particular, $S_0=\mu:=|A|\alpha\le d\alpha$.
At step $i$, only vertices in $A\cap\NR(v_i)$ can contribute to $S_i-S_{i-1}$. To control this increment and its conditional expectation under $\Qq$, we define $Y_i=\sum_{a\in A\cap\NR(v_i)}p_{i-1}(a)$ and $b_i=\E_{\Qq}[S_i-S_{i-1}\mid\F_{i-1}]$
for every $1\le i\le K$. Both $Y_i$ and $b_i$ are $\F_{i-1}$-measurable.
We next establish the following bounds for the weight processes under $\Qq$.

\begin{lemma}\label{lem:tilted-drift}
For every $a\in A$, the process $(p_i(a))_{i=0}^K$ is a nonnegative submartingale under $\Qq$. For every $1\le i\le K$, the bounds
\begin{equation}\label{eq:drift-bound}
0\le b_i\le\1_{\{v_i\in A\}}\frac{q_iY_i}{1-q_i}
\end{equation}
and
\begin{equation}\label{eq:increment-bound}
-Y_i\le S_i-S_{i-1}\le\frac{q_iY_i}{1-q_i}
\end{equation}
hold $\Qq$-almost surely.
\end{lemma}

\begin{proof}
Fix $1\le i\le K$ and $a\in A$. If $q_i=0$ or $a\notin N_R(v_i)$, then $p_i(a)=p_{i-1}(a)$. If $v_i\notin A$ or $p_{i-1}(v)=h$, then $p_i(v)=p_{i-1}(v)$, and~\eqref{eq:tilted-transition} together with~\eqref{eq:original-martingale} shows
$\E_{\Qq}[p_i(a)-p_{i-1}(a)\mid\F_{i-1}]=\E_{\Pp}[p_i(a)-p_{i-1}(a)\mid\F_{i-1}]=0$.

It remains to consider $v_i\in A$, $q_i>0$, and $x=p_{i-1}(v)\in(0,h)$. Fix $a\in A\cap\NR(v_i)$ with $y=p_{i-1}(a)<h$, and write $x'=p_i(v)$, $y'=p_i(a)$, and $q=q_i$. For $z\in\{x,y\}$, let $c_z$ and $\theta_z$ be the values of $c$ and $\theta$ in Algorithm~\ref{alg:random-independent-set} corresponding to the current weight $z$.
By Algorithm~\ref{alg:random-independent-set}, conditional on $\F_{i-1}$ and $A_i=0$, we have $x'=c_x$ and $y'=c_y$. 
Conditional on $\F_{i-1}$ and $A_i$, the variables $x'$ and $y'$ are independent. Their conditional means are $c_x,c_y$ when $A_i=0$, and $h\theta_x,h\theta_y$ when $A_i=1$. Therefore, $\Cov_{\Pp}(x',y'\mid\F_{i-1})=q(1-q)(c_x-h\theta_x)(c_y-h\theta_y)$. Since $c_z-h\theta_z=(c_z-z)/q$ and $0\le c_z-z\le qz/(1-q)$ for $z\in\{x,y\}$, we obtain
\begin{equation}\label{eq:covariance}
0\le\Cov_{\Pp}(x',y'\mid\F_{i-1})\le\frac{qxy}{1-q}.
\end{equation}

By~\eqref{eq:tilted-transition} and~\eqref{eq:original-martingale},
\[
\begin{aligned}
\E_{\Qq}[y'-y\mid\F_{i-1}]
&=\frac{\E_{\Pp}[x'y'\mid\F_{i-1}]}{x}-y
=\frac{\Cov_{\Pp}(x',y'\mid\F_{i-1})}{x}.
\end{aligned}
\]
Together with~\eqref{eq:covariance}, this gives
$
0\le\E_{\Qq}[y'-y\mid\F_{i-1}]
\le\frac{qy}{1-q}.
$
Thus $(p_i(a))_{i=0}^K$ is a nonnegative submartingale under $\Qq$ for every $a\in A$. Summing the above bounds over $a\in A$, we obtain~\eqref{eq:drift-bound}.

Finally, by Algorithm~\ref{alg:random-independent-set}, we have $0\le p_i(a)\le p_{i-1}(a)/(1-q_i)$ for every $a\in A\cap\NR(v_i)$. Since all other weights in $A$ remain unchanged, summing over $a\in A\cap\NR(v_i)$ proves~\eqref{eq:increment-bound}.
\end{proof}

\section{\texorpdfstring{$K_{t,t}$}{K(t,t)}-free graphs: Proof of Theorem~\ref{thm:main}}\label{sec:fractional-proof}
In this section, we prove Theorem~\ref{thm:main}. We first establish estimates for a stopped weight process using Lemma~\ref{lem:weighted} and martingale arguments. We then use these estimates to show that the probability of a vertex's weight reaching the upper threshold is small, which yields the desired lower bound on its selection probability.
\subsection{Running maxima and the stopped process}\label{subsec:stopping}
For each $a\in A$, the \emph{running maximum} of its weight up to step $i$ is $\max_{0\le r\le i}p_r(a)$. To control the total variation of the weights in $A$, we consider the sum of their running maxima:
$M_i=\sum_{a\in A}\max_{0\le r\le i}p_r(a).$
For every $a\in A$ and $0\le j\le K$,
\begin{equation}\label{eq:variation}
\sum_{i=1}^j|p_i(a)-p_{i-1}(a)|
\le2\max_{0\le r\le j}p_r(a).
\end{equation}
For each $a\in A$, we have
$
\max_{0\le r\le i}p_r(a)-\max_{0\le r\le i-1}p_r(a)
=\left(p_i(a)-\max_{0\le r\le i-1}p_r(a)\right)_+
\le(p_i(a)-p_{i-1}(a))_+.
$
By~\eqref{eq:weight-growth}, 
\begin{equation}\label{eq:max-overshoot}
M_i\le M_{i-1}+\frac{q_iY_i}{1-q_i}
\le\frac{M_{i-1}}{1-h},
\end{equation}
where the last inequality follows from $Y_i\le S_{i-1}\le M_{i-1}$ and $q_i\le h$.

Fix $R>\mu$ and $\delta>0$, and let
$
T=\min\bigl(\{i\in\{1,\ldots,K\}:M_i>R\}\cup\{K\}\bigr).
$
Thus $T$ is a stopping time. Moreover, $M_{T-1}\le R$, so~\eqref{eq:max-overshoot} implies
\begin{equation}\label{eq:stopped-max}
M_T\le\frac{R}{1-h}\le2R.
\end{equation}
We call step $i$ \emph{$\delta$-light} if $Y_i\le\delta$, and \emph{$\delta$-heavy} otherwise.

The following maximal inequality for nonnegative submartingales is classical; see~\cite[Section~14.6]{Williams}.

\begin{lemma}[Doob's maximal inequality \cite{Williams}]\label{lem:doob}
Let $(U_i)_{i=0}^K$ be a nonnegative submartingale with respect to $(\F_i)_{i=0}^K$ under $\Qq$. Then, for every $u>0$,
\[
u\,\Qq\left(\max_{0\le i\le K}U_i\ge u\right)
\le\E_{\Qq}\left[U_K\1_{\{\max_{0\le i\le K}U_i\ge u\}}\right]
\le\E_{\Qq}U_K.
\]
\end{lemma}
We use the following standard orthogonality identity for martingale differences; see~\cite[Section~12.1]{Williams}.

\begin{lemma}[orthogonality identity \cite{Williams}]\label{lem:orthogonal-differences}
Let $D_1,\ldots,D_K$ be square-integrable random variables such that $D_i$ is $\F_i$-measurable and $\E_{\Qq}[D_i\mid\F_{i-1}]=0$ for every $1\le i\le K$. Then $\E_{\Qq}[\sum_{i=1}^K D_i]=0$ and
\[
\E_{\Qq}\left(\sum_{i=1}^K D_i\right)^2=\sum_{i=1}^K\E_{\Qq}D_i^2.
\]
\end{lemma}

We next bound the total weight of $A$ and the expected sum of its running maxima up to the stopping time $T$.

\begin{proposition}\label{prop:stopped}
Suppose that $G$ is $K_{s+1,s+1}$-free, where $s\ge2$, and that
\begin{equation}\label{eq:threshold-stopped}
\delta\ge2s^{1/s}h^{1/s}(2R)^{1-1/s}.
\end{equation}
Let $\eta=\frac{2\delta R}{1-h}+\frac{8(s-1)hR^2}{\delta(1-h)}$. Then the following hold.
\begin{enumerate}
\item[(i)] There exists a random variable $Z$ with $\E_{\Qq}Z=0$ and $\E_{\Qq}Z^2\le2\delta\E_{\Qq}M_T\le4\delta R$ such that
\begin{equation}\label{eq:stopped-upper}
S_T\le\mu+Z+\eta.
\end{equation}

\item[(ii)] We have
\begin{equation}\label{eq:expected-max}
\E_{\Qq}M_T\le\mu+(\mu+\eta)\log(h/\alpha).
\end{equation}
\end{enumerate}
\end{proposition}

\begin{proof}
We begin by estimating the number of heavy steps before the process stops. For each $a\in A$, let $r_a=\max_{0\le i\le T}p_i(a)$. By~\eqref{eq:stopped-max}, we have $0\le r_a\le h$ and $\sum_{a\in A}r_a=M_T\le2R$. If $i\le T$ is a $\delta$-heavy step, then $\sum_{a\in A\cap N(v_i)}r_a\ge Y_i>\delta$. By Lemma~\ref{lem:weighted} and~\eqref{eq:threshold-stopped}, there are at most $2(s-1)M_T/\delta\le4(s-1)R/\delta$ such steps. 

By~\eqref{eq:increment-bound}, each $\delta$-heavy step $i\le T$ satisfies $(S_i-S_{i-1})_+\le q_iY_i/(1-q_i)\le2hR/(1-h)$. Hence
\begin{equation}\label{eta1}
\sum_{\substack{i\le T\\Y_i>\delta}}(S_i-S_{i-1})_+
\le\frac{2(s-1)M_T}{\delta}\cdot\frac{hM_T}{1-h}
=\frac{2(s-1)hM_T^2}{\delta(1-h)}
\le\frac{8(s-1)hR^2}{\delta(1-h)}.
\end{equation}
For $\delta$-light steps,~\eqref{eq:drift-bound} shows
\begin{equation}\label{eta2}
\sum_{\substack{i\le T\\Y_i\le\delta}}b_i\le\frac{\delta}{1-h}\sum_{\substack{i\le T\\v_i\in A}}q_i\le\frac{\delta M_T}{1-h}\le\frac{2\delta R}{1-h},
\end{equation}
where we used $q_i\le r_{v_i}$ whenever $i\le T$ and $v_i\in A$.

Define
$
Z=\sum_{i=1}^K\1_{\{i\le T,\,Y_i\le\delta\}}(S_i-S_{i-1}-b_i).
$
Since $T$ is a stopping time and $Y_i$ is $\F_{i-1}$-measurable, the event $\{i\le T,\,Y_i\le\delta\}$ belongs to $\F_{i-1}$. Hence $Z$ is a sum of martingale differences under $\Qq$, and $\E_{\Qq}Z=0$.
For $\delta$-light steps,~\eqref{eq:increment-bound} and $q_i\le h<1/2$ give $|S_i-S_{i-1}|\le\delta$. 
By Lemma \ref{lem:orthogonal-differences} and \eqref{eq:variation}, we obtain
\[
\begin{aligned}
\E_{\Qq}Z^2
\le\E_{\Qq}\sum_{\substack{i\le T\\Y_i\le\delta}}(S_i-S_{i-1})^2
\le\delta\E_{\Qq}\sum_{i\le T}\sum_{a\in A}|p_i(a)-p_{i-1}(a)|
\le2\delta\E_{\Qq}M_T\le4\delta R.
\end{aligned}
\]
By the definition of $Z$, we have $\sum_{i\le T,\,Y_i\le\delta}(S_i-S_{i-1})=Z+\sum_{i\le T,\,Y_i\le\delta}b_i$. Moreover, $\sum_{i\le T,\,Y_i>\delta}(S_i-S_{i-1})\le\sum_{i\le T,\,Y_i>\delta}(S_i-S_{i-1})_+$. 
By \eqref{eta1} and \eqref{eta2}, we have
\[
\begin{aligned}
S_T
&=\mu+\sum_{\substack{i\le T\\Y_i\le\delta}}(S_i-S_{i-1})
+\sum_{\substack{i\le T\\Y_i>\delta}}(S_i-S_{i-1})\\
&\le\mu+Z+\sum_{\substack{i\le T\\Y_i\le\delta}}b_i
+\sum_{\substack{i\le T\\Y_i>\delta}}(S_i-S_{i-1})_+\\
&\le\mu+Z+\eta,
\end{aligned}
\]
which proves \eqref{eq:stopped-upper}.

Finally, fix $a\in A$. By Lemma~\ref{lem:tilted-drift}, the stopped process $(p_{\min\{i,T\}}(a))_{i=0}^K$ is a nonnegative submartingale under $\Qq$. By Lemma~\ref{lem:doob}, for every $u\in[\alpha,h]$,
$
\Qq\left(\max_{0\le i\le T}p_i(a)\ge u\right)\le\frac{\E_{\Qq}p_T(a)}{u}.
$
Since $p_0(a)=\alpha$ and $p_i(a)\le h$, integrating over $u\in[\alpha,h]$ shows
$
\E_{\Qq}\max_{0\le i\le T}p_i(a)\le\alpha+\E_{\Qq}p_T(a)\log(h/\alpha).
$
Summing over $a\in A$ and taking expectations in~\eqref{eq:stopped-upper}, using $\E_{\Qq}Z=0$, we obtain~\eqref{eq:expected-max}.
\end{proof}

%\section{Proof of the fractional-coloring theorem}\label{sec:fractional-proof}
\subsection{Completion of the proof}

We first establish a lower bound on $S_K$ whenever $\mathcal B_k$ occurs.

\begin{lemma}\label{lem:cap}
On the event $\mathcal B_k$, we have $S_K\ge(1-h)\log(h/\alpha).$
\end{lemma}

\begin{proof}
The weight of $v$ changes only when a vertex of $A$ is processed. Since $q_i\le h$, we have $(1-q_i)^{-1}\le\exp(q_i/(1-h))$. Moreover, $q_i\le p_K(v_i)$ for every $v_i\in A$, since the weight of a vertex remains unchanged after it is processed. By~\eqref{eq:weight-growth},
\[
p_K(v)\le\alpha\prod_{v_i\in A}\frac{1}{1-q_i}
\le\alpha\exp\left(\frac{1}{1-h}\sum_{v_i\in A}q_i\right)
\le\alpha\exp\left(\frac{S_K}{1-h}\right).
\]
Since $p_K(v)=h$ on $\mathcal B_k$, the conclusion follows.
\end{proof}

Now we are ready to prove Theorem~\ref{thm:main}.

\begin{proof}[Proof of Theorem~\ref{thm:main}]

First suppose that $t\ge3$. We choose the parameters so that the probability of stopping early is small and the total weight of the left-neighbors is unlikely to increase substantially. 
Choose $h=(\log d)^{-16(t-1)}$, $R=(\log d)^4$, $\delta=(\log d)^{-8}$ and $\alpha=(\log d-(16t-13)\log\log d)/d$. For sufficiently large $d$ depending only on $t$, we have $0<\alpha<h<1/2$ and $\mu\le d\alpha<\log d<R$. 
To apply Proposition~\ref{prop:stopped}, take $s=t-1$. Since $h^{1/(t-1)}=(\log d)^{-16}$, the right-hand side of~\eqref{eq:threshold-stopped} equals $2^{2-1/(t-1)}(t-1)^{1/(t-1)}(\log d)^{-12-4/(t-1)}=o(\delta)$. Thus the required condition holds for sufficiently large $d$. Moreover, by the definition of $\eta$,
$
\eta=\frac{2(\log d)^{-4}+8(t-2)(\log d)^{-16(t-2)}}{1-(\log d)^{-16(t-1)}}=O_t((\log d)^{-4}).
$
In particular, $\eta\le1$ for sufficiently large $d$.

Since $\log(h/\alpha)\le\log d$, by Proposition~\ref{prop:stopped}, we have $\E_{\Qq}M_T\le\log d+(\log d+1)\log d\le2(\log d)^2$. 
By Markov's inequality and the bound $\E_{\Qq}M_T\le2(\log d)^2$, we obtain
\begin{equation}\label{eq:stop-prob}
\Qq(M_T>R)\le\frac{\E_{\Qq}M_T}{R}\le\frac{2}{(\log d)^2}.
\end{equation}
Moreover, Proposition~\ref{prop:stopped} implies
$
\E_{\Qq}Z^2\le4\delta R=4(\log d)^{-4}.
$

By our choices of $h$ and $\alpha$, we have $\log d-(16t-15)\log\log d\le\log(h/\alpha)\le\log d$. Since $\mu\le\log d-(16t-13)\log\log d$, it follows that
$
(1-h)\log(h/\alpha)-\mu-\eta
\ge 2\log\log d-h\log d-\eta
\ge\log\log d
$
for sufficiently large $d$.

On the event $M_T\le R$, we have $T=K$. By Lemma~\ref{lem:cap} and~\eqref{eq:stopped-upper}, the event $\mathcal B_k$ then implies $Z\ge\log\log d$. 
Since $\mathcal B_k\subseteq\{M_T>R\}\cup\{Z\ge\log\log d\}$, combining~\eqref{eq:stop-prob} with Chebyshev's inequality and $\E_{\Qq}Z^2\le4(\log d)^{-4}$ gives
\[
\Qq(\mathcal B_k)
\le\frac{2}{(\log d)^2}
+\frac{4}{(\log d)^4(\log\log d)^2}
\le\frac{4}{(\log d)^2}.
\]
Thus, by~\eqref{eq:selection-Q}, every vertex $v$ with $A\ne\varnothing$ satisfies
\[
\Pp(v\in I)\ge\alpha\left(1-\frac{4}{(\log d)^2}\right)
\ge\frac{\log d-(16t-12)\log\log d}{d}.
\]
The same bound holds when $A=\varnothing$, since $\Pp(v\in I)=\alpha$. By the probability characterization of fractional coloring, this proves~\eqref{eq:main} for $t\ge3$.

For $t=2$, use the same parameters. Since $G$ is $C_4$-free, every vertex other than $v$ has at most one neighbor in $A$, so $Y_i\le h<\delta$. The proof of Proposition~\ref{prop:stopped} therefore applies without the heavy-step estimate, with $\eta=2\delta R/(1-h)$. The preceding argument gives $\Pp(v\in I)\ge(\log d-20\log\log d)/d$ for every vertex $v$, proving~\eqref{eq:main}.
%For $t=2$, choose the parameters as above. Since $G$ is $C_4$-free, every vertex other than $v$ has at most one neighbor in $A$. Hence $Y_i\le h<\delta$ for every $1\le i\le K$, so all steps are $\delta$-light. By~\eqref{eq:drift-bound}, we have $\sum_{i\le T}b_i\le\delta M_T/(1-h)\le\eta$, where $\eta=2\delta R/(1-h)=O((\log d)^{-4})$.
%Let $Z=\sum_{i=1}^K\1_{\{i\le T\}}(S_i-S_{i-1}-b_i)$. Since $\{i\le T\}\in\F_{i-1}$, the summands are martingale differences under $\Qq$, so $\E_{\Qq}Z=0$. By~\eqref{eq:increment-bound}, we have $|S_i-S_{i-1}|\le\delta$. Using the total variation estimate, we obtain $\E_{\Qq}Z^2\le2\delta\E_{\Qq}M_T\le4\delta R$. Consequently, $S_T\le\mu+Z+\eta$ and $\E_{\Qq}S_T\le\mu+\eta$. Applying Lemma~\ref{lem:doob} to the stopped weight processes gives $\E_{\Qq}M_T\le\mu+(\mu+\eta)\log(h/\alpha)\le2(\log d)^2$.
%It follows that $\Qq(M_T>R)\le2/(\log d)^2$. Moreover, with $B'=(1-h)\log(h/\alpha)$, we have $B'-\mu-\eta\ge\log\log d$ for sufficiently large $d$. By Lemma~\ref{lem:cap}, the event $\mathcal B_k$ therefore implies either $M_T>R$ or $Z\ge\log\log d$. Chebyshev's inequality yields $\Qq(\mathcal B_k)\le4/(\log d)^2$. Hence $\Pp(v\in I)\ge(\log d-20\log\log d)/d$, which proves~\eqref{eq:main} for $t=2$.
\end{proof}

\begin{remark}
The independent-set sampling procedure also yields an efficient randomized construction of a fractional coloring. As in the proof of~\cite[Theorem~1.3]{ADN}, independently sampling sufficiently many independent sets and applying a Chernoff bound gives, with high probability, a fractional coloring of total weight at most $(1+\varepsilon)d/\log d$, supported on $O_{t,\varepsilon}(d\log n/\log d)$ sampled independent sets.
\end{remark}

We now deduce Theorem~\ref{thm:bipartite} from Theorem~\ref{thm:main} using Lemma \ref{lem:bipartite-transfer}.

\begin{proof}[Proof of Theorem~\ref{thm:bipartite} assuming Theorem~\ref{thm:main}]
By Lemma~\ref{lem:degenerate-reduction}, $G$ contains a $d$-degenerate induced subgraph $H$ with $\delta(H)\ge d$. Since $H$ is also $K_{t,t}$-free, Theorem~\ref{thm:main} shows $\chif(H)\le d/(\log d-16t\log\log d)$. As $\overline d(H)\ge d$, we may apply Lemma~\ref{lem:bipartite-transfer} to obtain an induced bipartite subgraph of $G$ with average degree at least
$
\frac{\overline d(H)}{\chif(H)}\ge\log d-16t\log\log d,
$
as required.
\end{proof}

\section{Triangle-free graphs: Proof of Theorem~\ref{thm:triangle-free}}\label{sec:triangle-free}
Our proof combines an entropy estimate with a stopping rule for the independent-set procedure in Section~\ref{sec:sampler}. The entropy estimate improves the bound obtained from the procedure of Martinsson and Steiner~\cite[Section~4]{MS} by repeated application. The stopping rule controls the weighted sum over a right-neighborhood. We then compare two bounds for the same weighted expression, one from entropy and one from a submartingale. This comparison gives the leading constant $2.85$. %The entropy and edge-product calculations are closely related to those of Dhawan~\cite[Lemma~2.3]{Dhawan}.

\medskip

\begin{algorithm}[H]
\caption{The Martinsson--Steiner procedure~\cite[Section~4]{MS}}\label{alg:tf-ms}
\KwIn{A graph $G$, a vertex ordering $v_1,\ldots,v_n$, and initial weights $x_v\ge0$.}
$I\gets\varnothing$ and $p_0(v)\gets x_v$ for every $v\in V(G)$\;
\For{$i=1,\ldots,n$}{
  $P_{v_i}\gets p_{i-1}(v_i)$\;
  Choose $A_i\in\{0,1\}$ with $\Pp(A_i=1\mid\F_{i-1})=1-e^{-P_{v_i}}$\;
  \If{$A_i=1$}{$I\gets I\cup\{v_i\}$\;}
  \For{$u\in\NR(v_i)$}{
    \eIf{$A_i=1$}{
      $p_i(u)\gets0$\;
    }{
      $p_i(u)\gets p_{i-1}(u)e^{P_{v_i}}$\;
    }
  }
  Leave all other weights unchanged\;
}
\end{algorithm}

\subsection{An entropy estimate}

For a nonnegative vertex weight function $w:V(G)\to[0,\infty)$, define
\[
\alpha_w(G)=\max\left\{\sum_{v\in I}w_v:I\subseteq V(G)\text{ is independent}\right\}.
\]
We also define $H(x)=x\log(1/x)$ for $x>0$, with $H(0)=0$. Let $\kappa=1/[2(\sqrt e-1)]$ and $\Psi(x)=2\kappa\log(1+x/(2\kappa))$ for $0\le x\le1$. Note that $\Psi(1)=\kappa$.

We begin with the following entropy estimate, which will play a key role in the proof.

\begin{lemma}\label{lem:tf-record-entropy}
Let $G$ be a triangle-free graph with a fixed vertex ordering, and let $Z>0$. Suppose that $w_v\ge0$, $0\le q_v<1$, and $z_v\ge Z$ are random variables indexed by $V(G)$, with $\E\alpha_w(G)\le1$ and $A:=\E\sum_vw_vq_v\le1$. Suppose also that, almost surely, if $q_v>0$, then 
\begin{equation}\label{eq:tf-record-assumptions}
\log(1/q_v)+\sum_{u\in\NL(v)}q_u\le z_v.
\end{equation}
Then $\E\sum_{q_v>0}w_vH(q_v)/z_v\le\Psi(A)+o_Z(1)$.
\end{lemma}

\begin{proof}
Define $b_0=1$ and $b_{j+1}=2b_j\log(1+1/(2b_j))$ for $j\ge0$. We first prove by induction on $j$ that
\begin{equation}\label{eq:tf-basic-entropy}
\E\sum_{x_v>0}w_vx_v\left(1-\frac{\sum_{u\in\NL(v)}x_u}{\log(1/x_v)}\right)\le b_j+o_Z(1)
\end{equation}
whenever the random vertex weights $w_v\ge0$ and $0\le x_v\le e^{-Z}$ satisfy $\E\alpha_w(G)\le1$ and $\E\sum_vw_vx_v\le M$, where $M>0$ is fixed.

For $j=0$, condition on $(w,x)$ and apply Algorithm~\ref{alg:tf-ms}. By construction, for every $u\in\NR(v_i)$, we have $\E[p_i(u)\mid\F_{i-1}]=p_{i-1}(u)$ and $\E[H(p_i(u))\mid\F_{i-1}]=H(p_{i-1}(u))-P_{v_i}p_{i-1}(u)$. Since $G$ is triangle-free, at most one endpoint of an edge changes weight at each step. Thus $p_i(u)p_i(v)$ is a martingale for every $uv\in E(G)$. It follows that $\E[P_v\mid w,x]=x_v$ and $\E[H(P_v)\mid w,x]=H(x_v)-x_v\sum_{u\in\NL(v)}x_u$.
For $x_v>0$, we have $P_v\in\{0\}\cup[x_v,\infty)$ and
$H(P_v)/\log(1/x_v)\le1-e^{-P_v}+P_v/\log(1/x_v)$. Consequently,
\[
\begin{aligned}
&\mathbb E\sum_{x_v>0}w_vx_v\left(1-\frac{\sum_{u\in N_L(v)}x_u}{\log(1/x_v)}\right)\\
&\quad=\mathbb E\sum_{x_v>0}\frac{w_vH(P_v)}{\log(1/x_v)}\\
&\quad\le\mathbb E\sum_vw_v(1-e^{-P_v})+\frac1Z\mathbb E\sum_vw_vP_v\\
&\quad=\mathbb E\sum_{v\in I}w_v+\frac1Z\mathbb E\sum_vw_vx_v\\
&\quad\le\mathbb E\alpha_w(G)+\frac MZ\le1+\frac MZ.
\end{aligned}
\]
This proves the case $j=0$.

Suppose that~\eqref{eq:tf-basic-entropy} holds with coefficient $b_j$ for some $j\ge0$. We first establish the following estimate, which will be used to prove the same inequality with coefficient $b_{j+1}$.
We first show that
\begin{equation}\label{eq:tf-entropy-improvement}
\E\sum_{q_v>0}\frac{w_vH(q_v)}{z_v}
\le 2b_j\log\left(1+\frac{A}{2b_j}\right)+o_Z(1).
\end{equation}

For each $v$ with $q_v>0$, let $t_v=z_v/\log(1/q_v)-1$. By~\eqref{eq:tf-record-assumptions}, we have $t_v\ge0$ and $\sum_{u\in\NL(v)}q_u\le t_v\log(1/q_v)$. For $t\ge0$, define $J(t)=\E\sum_{q_v>0}w_vq_v(t-t_v)_+$.

Fix $t>0$, and let $x_v=q_v/t$ if $q_v>0$ and $t_v<t$, and $x_v=0$ otherwise. For every $v$ with $x_v>0$, we have $\log(1/q_v)\ge Z/(1+t)$ and $\sum_{u\in\NL(v)}x_u\le(t_v/t)\log(1/q_v)$. Moreover, $\E\sum_vw_vx_v\le1/t$. Since $\log(1/x_v)=\log(1/q_v)+\log t$, replacing $\log(1/x_v)$ by $\log(1/q_v)$ in~\eqref{eq:tf-basic-entropy} introduces an error of $O_t(1/Z)$. Applying the induction hypothesis with $Z$ replaced by $Z/(1+t)+\log t$, we obtain $J(t)\le b_jt^2+o_Z(1)$.

Since $J(0)=0$ and $|J(t)-J(s)|\le A|t-s|$, this bound is uniform on bounded intervals. Taking $s=\min\{t,A/(2b_j)\}$ gives $J(t)\le b_js^2+A(t-s)+o_Z(1)$ uniformly for $t\ge0$. Using the identity $(1+t_v)^{-1}=\int_0^\infty 2(t-t_v)_+/(1+t)^3\,dt$, we conclude that
\[
\E\sum_{q_v>0}\frac{w_vH(q_v)}{z_v}
=\int_0^\infty\frac{2J(t)}{(1+t)^3}\,dt
\le 2b_j\log\left(1+\frac{A}{2b_j}\right)+o_Z(1).
\]
Here we used $J(t)\le b_jt^2+o_Z(1)$ for $t\le A/(2b_j)$ and $J(t)\le At-A^2/(4b_j)+o_Z(1)$ otherwise. This proves~\eqref{eq:tf-entropy-improvement}.

It remains to prove~\eqref{eq:tf-basic-entropy} with coefficient $b_{j+1}$. Run Algorithm~\ref{alg:tf-ms} with initial weights $x_v$. Let $\delta=Z^{-1/2}$ and $q_v=P_v\mathbf1_{\{0<P_v\le\delta\}}$. Thus, we have $\sum_{u\in\NL(v)}P_u=\log(P_v/x_v)$. If $q_v>0$, then $\log(1/q_v)+\sum_{u\in\NL(v)}q_u\le\log(1/x_v)$.
Taking $z_v=\log(1/x_v)$ when $x_v>0$ and $z_v=Z$ otherwise, we see that~\eqref{eq:tf-record-assumptions} is satisfied.

Since $p\le(1+\delta)(1-e^{-p})$ for $0\le p\le\delta$, we have $\E\sum_vw_vq_v\le(1+\delta)\E\sum_{v\in I}w_v\le1+\delta$. Applying~\eqref{eq:tf-entropy-improvement} with weights $w_v/(1+\delta)$, we obtain
\[
\E\sum_{0<P_v\le\delta}\frac{w_vH(P_v)}{\log(1/x_v)}
\le(1+\delta)b_{j+1}+o_Z(1)
=b_{j+1}+o_Z(1).
\]
The contribution from vertices with $P_v>\delta$ is at most $M\log(1/\delta)/Z=o_Z(1)$, since $H(P_v)\le P_v\log(1/\delta)$ and $\E[P_v\mid w,x]=x_v$. The entropy identity therefore gives~\eqref{eq:tf-basic-entropy} with coefficient $b_{j+1}$, completing the induction.

Finally, the sequence $(b_j)$ decreases to $\kappa$. Indeed, the function $b\mapsto2b\log(1+1/(2b))$ has $\kappa$ as its unique positive fixed point and is strictly smaller than $b$ for $b>\kappa$. For each fixed $j$, the preceding argument gives the conclusion of the lemma with $2b_j\log(1+A/(2b_j))$ in place of $\Psi(A)$. As $j\to\infty$, these bounds converge uniformly to $\Psi(A)$ for $0\le A\le1$. The conclusion follows by choosing $j$ sufficiently large and then letting $Z\to\infty$.
\end{proof}

\subsection{A stopping rule}

By linear programming duality, $\chif(G)=\max\{\sum_vw_v:w_v\ge0,\ \alpha_w(G)\le1\}$; see~\cite[Section~3.1]{SU}. We may assume that $G$ is nonempty. Choose a maximizing weight function $w$, and write $W=\sum_vw_v=\chif(G)$. Then $\alpha_w(G)=1$. Deleting vertices of weight zero, we may assume that $w_v>0$ for every $v$.

Fix a $d$-degeneracy ordering. We use the same initial weight at every vertex, choosing $\alpha$ so that $d\alpha=(5/8)\log(1/\alpha)$. Write $z=\log(1/\alpha)$ and $h=z^{-2}$. Here $z$ is the initial logarithmic weight, while $h$ bounds every weight during the procedure. For sufficiently large $d$, we have $\alpha<h<1/2$ and $z\sim\log d$. Theorem~\ref{thm:MS-degenerate} gives $\alpha W=O(1)$.

We modify Algorithm~\ref{alg:random-independent-set} by stopping vertices according to $U_i(v)$ and $S_i(v)$, as described in Algorithm~\ref{alg:tf-stopping}. For the weight process, define $U_i(v)=(zw_v)^{-1}\sum_{u\in\NR(v)}w_up_i(u)$. When $p_i(v)>0$, define $S_i(v)=z^{-1}\log(1/p_i(v))$, and otherwise set $S_i(v)=0$. The quantity $U_i(v)$ measures the weighted sum over the later neighbors, normalized so that it controls the entropy lost when $v$ is processed. The quantity $S_i(v)$ starts at $1$ and decreases as the positive weight of $v$ increases.

\begin{algorithm}[tbp]
\caption{The independent set procedure with a stopping rule}\label{alg:tf-stopping}
\KwIn{$G$, an ordering $v_1,\ldots,v_n$, vertex weights $w_v>0$, and parameters $0<\alpha<h<1/2$.}
\KwOut{An independent set $I$, selection probabilities $(q_v)$, and weights $(p_i(v))$.}
$I\gets\varnothing$ and $p_0(v)\gets\alpha$ for every $v\in V(G)$\;
$\mathcal A_0\gets\{v\in V(G):U_0(v)\le7/4\}$\;
\For{$i=1,\ldots,n$}{
  $p_i(v)\gets p_{i-1}(v)$ for every $v\in V(G)$\;
  $q_{v_i}\gets p_{i-1}(v_i)\mathbf1_{\{v_i\in\mathcal A_{i-1}\}}$\;
  \If{$q_{v_i}>0$}{
    Choose $B_i\in\{0,1\}$ with $\Pp(B_i=1\mid\F_{i-1})=q_{v_i}$\;
    \If{$B_i=1$}{$I\gets I\cup\{v_i\}$\;}
    \ForEach{$v\in\mathcal A_{i-1}\cap\NR(v_i)$}{
      $x\gets p_{i-1}(v)$, $b\gets\min\{x/(1-q_{v_i}),h\}$\;
      $\theta\gets[x-(1-q_{v_i})b]/(q_{v_i}h)$\;
      \eIf{$B_i=0$}{
        $p_i(v)\gets b$\;
      }{
        Independently set $p_i(v)\gets h$ with probability $\theta$, and $p_i(v)\gets0$ otherwise\;
      }
    }
  }
  $\mathcal A_i\gets\{v\in\mathcal A_{i-1}\setminus\{v_i\}:0<p_i(v)<h,\ U_i(v)\le1+3S_i(v)/4\}$\;
}
\end{algorithm}

For $v=v_j$, let $q_v$ be its conditional selection probability at step $j$, given $\F_{j-1}$. Let $\tau_v$ be the first time $i<j$ at which $v$ is stopped, taking $\tau_v=j-1$ if there is no such time. Write $r_v=p_{\tau_v}(v)$, $U_v=U_{\tau_v}(v)$, and $S_v=S_{\tau_v}(v)$. The weight of $v$ remains unchanged after $\tau_v$, so $r_v=p_n(v)$. If $v$ is active when processed, then $q_v=r_v$; otherwise $q_v=0$.

As in Section~\ref{sec:sampler}, the selected set $I$ is independent, every weight is a martingale, and all weights lie in $\{0\}\cup[\alpha,h]$. In particular, $\E r_v=\alpha$ and $0\le S_v\le1$. Triangle-freeness also implies that the product of the weights at two adjacent vertices is a martingale, since at most one of the two weights changes in any step. Finally,
\begin{equation}\label{eq:tf-selection-weight}
\E\sum_vw_vq_v=\E\sum_{v\in I}w_v\le1.
\end{equation}

\subsection{Two bounds for the stopped process}
The next two lemmas estimate the same weighted expression. Its entropy term is controlled by Lemma~\ref{lem:tf-record-entropy}. Its neighborhood term allows us to account for vertices removed by the stopping rule. 

\begin{lemma}\label{lem:tf-stopped-estimate}
In Algorithm \ref{alg:tf-stopping}, we have
\[
\E\sum_vw_vq_v[(7/2-U_v)S_v+3U_v/4]\ge\frac{99}{64}\alpha W-o(1).
\]
\end{lemma}

\begin{proof}
We first estimate the entropy loss during the sampling procedure. At step $i$, write $q=q_{v_i}$ and $x=p_{i-1}(u)$ for a vertex $u$ whose weight is updated. Since $\E[p_i(u)\mid\F_{i-1}]=x$ and $p_i(u)\le x/(1-q)$, we have $\E[H(p_i(u))\mid\F_{i-1}]\ge H(x)+x\log(1-q)\ge H(x)-qx/(1-h)$.

Summing over all steps and vertices, and using $H(\alpha)=\alpha z$ and $H(r_v)=zr_vS_v$, we obtain
$\E\sum_vw_vr_vS_v\ge\alpha W-(1-h)^{-1}\E\sum_vw_vq_vU_v$.
Here the sum over updated vertices is bounded by the sum over all later neighbors. Since $\E r_v=\alpha$, $U_v\le7/4$ whenever $q_v>0$, and $\E\sum_vw_vq_v\le1$, it follows that
\begin{equation}\label{eq:tf-total-entropy-loss}
\E\sum_vw_v(r_vS_v+q_vU_v-r_v)
\ge-\frac{7h}{4(1-h)}=-O(h).
\end{equation}

For the rest of the proof, define $V_i(v)=z^{-1}\sum_{u\in\NL(v)}p_i(u)$ and $V_v=V_{\tau_v}(v)$. Both $U_i(v)$ and $V_i(v)$ are martingales. For each vertex $v$, define
$X_i(v)=p_i(v)(7/4-U_i(v))(2-V_i(v))$.
We claim that $(X_i(v))_{i=0}^n$ is a submartingale. Fix $i$. If $v_i\in N(v)$, then triangle-freeness implies that $U_i(v)$ and $V_i(v)$ remain unchanged. Since $p_i(v)$ is a martingale, we have $\E[X_i(v)\mid\F_{i-1}]=X_{i-1}(v)$. Otherwise, $p_i(v)=p_{i-1}(v)$, and
$\E[X_i(v)\mid\F_{i-1}]-X_{i-1}(v)=p_{i-1}(v)\Cov(U_i(v),V_i(v)\mid\F_{i-1})\ge0$
by~\eqref{eq:covariance}. Thus $(X_i(v))_{i=0}^n$ is a submartingale. Applying the optional stopping theorem at $\tau_v$ (see~\cite[Sections~10.9--10.10]{Williams}), we obtain
\begin{equation}\label{eq:tf-stopped-potential}
\E\bigl[r_v(7/4-U_v)(2-V_v)\bigr]
\ge\alpha(7/4-U_0(v))(2-V_0(v)).
\end{equation}

We next relate the preceding submartingale estimate to the terms involving $q_v$. For each $v\in\mathcal A_0$, we show that
\begin{equation}\label{eq:tf-pointwise-potential}
\begin{aligned}
q_v[(7/2-U_v)S_v+3U_v/4]
&\ge\frac74(r_vS_v+q_vU_v-r_v)+r_v(7/4-U_v)(2-V_v)\\
&\quad-O\left(r_v\left[h(1+V_v)+\frac{\log(1/h)}{z}\right]\right).
\end{aligned}
\end{equation}

The inequality is immediate if $r_v=0$. Suppose that $r_v>0$. By the weight growth bound, $V_v\ge(1-h)(1-S_v)$. Moreover, the stopping rule and the bound on each update give $U_v\le7/[4(1-h)]$. Since $0\le S_v\le1$, these inequalities imply
$(7/4-U_v)(V_v+S_v-1)\ge-O(h(1+V_v))$.

If $q_v=r_v$, the difference between the left-hand side and the two main terms on the right-hand side of~\eqref{eq:tf-pointwise-potential} is exactly $r_v(7/4-U_v)(V_v+S_v-1)$. If $q_v=0$, this difference has the additional term $r_v[U_v(1+S_v)-(7/2)S_v]$. When $0<r_v<h$, the stopping rule gives $U_v>1+3S_v/4$, so this additional term is nonnegative. When $r_v=h$, it is at least $-(7/2)r_v\log(1/h)/z$. This proves~\eqref{eq:tf-pointwise-potential}.

For $v\notin\mathcal A_0$, we have $q_v=0$, $r_v=\alpha$, and $S_v=1$, so these vertices contribute nothing to the left-hand side or to $r_vS_v+q_vU_v-r_v$. The total expected error in~\eqref{eq:tf-pointwise-potential} is $o(1)$, since $\E\sum_vw_vr_v=\alpha W=O(1)$ and, by the edge-product martingale and optional stopping,
$\E\sum_vw_vr_vV_v=\alpha\sum_vw_vV_0(v)\le(5/8)\alpha W=O(1)$.

Summing~\eqref{eq:tf-pointwise-potential} over $v\in\mathcal A_0$ and applying~\eqref{eq:tf-total-entropy-loss} and~\eqref{eq:tf-stopped-potential}, we obtain
\[
\E\sum_vw_vq_v[(7/2-U_v)S_v+3U_v/4]
\ge\alpha\sum_{v\in\mathcal A_0}w_v(7/4-U_0(v))(2-V_0(v))-o(1).
\]
Finally, $V_0(v)\le5/8$ for every vertex, and $\sum_vw_vU_0(v)=\sum_vw_vV_0(v)\le5W/8$. Since $U_0(v)>7/4$ outside $\mathcal A_0$, the last sum is at least
$\frac{11}{8}\alpha\sum_{v\in\mathcal A_0}w_v(7/4-U_0(v))\ge\frac{11}{8}\alpha(7/4-5/8)W=\frac{99}{64}\alpha W$.
This completes the proof.
\end{proof}

\begin{lemma}\label{lem:tf-entropy-upper}
We have
\[
\E\sum_vw_vq_v[(7/2-U_v)S_v+3U_v/4]\le\frac{33}{16}\kappa+\frac{233}{200}+o(1).
\]
\end{lemma}

\begin{proof}
For $q_v>0$, the construction gives $q_v=\alpha\prod_{u\in\NL(v)}(1-q_u)^{-1}$ and $S_v=\log(1/q_v)/z$. Hence, $\log(1/q_v)+\sum_{u\in\NL(v)}q_u\le\log(1/\alpha)=z$. Thus Lemma~\ref{lem:tf-record-entropy}, applied with $z_v=Z=z$ and using~\eqref{eq:tf-selection-weight}, gives $\E\sum_vw_vq_vS_v\le\Psi(1)+o(1)=\kappa+o(1)$.

For $q_v>0$, let $\theta_v=U_v/(1+3S_v/4)$, and let $\theta_v=0$ otherwise. Since $U_v=\theta_v(1+3S_v/4)$ whenever $q_v>0$, we have
\[
(7/2-U_v)S_v+3U_v/4
\le\frac{33}{16}S_v+\frac{23}{16}\left[(1-\theta_v)S_v+\frac{52}{69}\theta_v\right].
\]
Indeed, the difference between the right and left sides equals $3\theta_v(S_v-2/3)^2/4\ge0$.

Let $A=\E\sum_vw_vq_v(1-\theta_v)$. Since $0\le\theta_v\le1$, the modified weights $w'_v=w_v(1-\theta_v)$ satisfy $\alpha_{w'}(G)\le\alpha_w(G)=1$. Applying Lemma~\ref{lem:tf-record-entropy} with weights $w'_v$ and $z_v=Z=z$, we obtain $\E\sum_vw_vq_v(1-\theta_v)S_v\le\Psi(A)+o(1)$. Moreover, by~\eqref{eq:tf-selection-weight}, we have $\E\sum_vw_vq_v\theta_v\le1-A$. Consequently,
\[
\begin{aligned}
&\E\sum_vw_vq_v\left[(1-\theta_v)S_v+\frac{52}{69}\theta_v\right]
\le\Psi(A)+\frac{52}{69}(1-A)+o(1)
\le\frac{466}{575}+o(1).
\end{aligned}
\]
For the last inequality, the maximum occurs at $A=17\kappa/26$, where $\Psi(A)+(52/69)(1-A)=52/69+2\kappa(\log(69/52)-17/69)<466/575$.

Combining these estimates yields
\[
\begin{aligned}
&\E\sum_vw_vq_v[(7/2-U_v)S_v+3U_v/4]
\le\frac{33}{16}\kappa+\frac{23}{16}\cdot\frac{466}{575}+o(1)
=\frac{33}{16}\kappa+\frac{233}{200}+o(1),
\end{aligned}
\]
as required.
\end{proof}

\subsection{Completion of the proof}

\begin{proof}[Proof of Theorem~\ref{thm:triangle-free}]
Combining Lemmas~\ref{lem:tf-stopped-estimate} and~\ref{lem:tf-entropy-upper} gives $(99/64)\alpha W\le33\kappa/16+233/200+o(1)$. Since $d\alpha=5z/8$, it follows that
\[
W\le\left[\frac{512}{495}\left(\frac{33}{16}\kappa+\frac{233}{200}\right)+o(1)\right]\frac{d}{z}.
\]
Since $\frac{512}{495}\left(\frac{33}{16}\kappa+\frac{233}{200}\right)<\frac{57}{20}$ and $z\sim\log d$, the result follows.
\end{proof}

We finish by deducing the induced bipartite subgraph bound.

\begin{proof}[Proof of Theorem~\ref{thm:triangle-free-bipartite} assuming Theorem~\ref{thm:triangle-free}]
By Lemma~\ref{lem:degenerate-reduction}, $G$ contains a $d$-degenerate induced subgraph $H$ with $\delta(H)\ge d$. Since $H$ is triangle-free, Theorem~\ref{thm:triangle-free} gives $\chif(H)\le(57/20+o(1))d/\log d$. Lemma~\ref{lem:bipartite-transfer} therefore yields an induced bipartite subgraph of $G$ with average degree at least $\overline d(H)/\chif(H)\ge(20/57-o(1))\log d$.
\end{proof}

\section{bounded edge-codegrees}\label{sec:edge}

%We now prove Theorem~\ref{thm:edge-codegree}, based on the ideas~\cite[proof of Theorem~2]{Glock}. 

%\subsection{A joint binomial estimate}
Our proof builds on the proof of Theorem~2 in Glock~\cite{Glock}.
We first estimate the probability that a random subset contains exactly $r$ elements from each of two given sets.

\begin{lemma}\label{lem:joint-binomial}
Let $d\ge2$ be an integer, and let $A$ and $B$ be sets with $|A|=|B|=d$. Write $b=|A\cap B|$, and let $R$ be a random subset of $A\cup B$ obtained by selecting each element independently with probability $p=1/d$. For an integer $1\le r\le d/2$, let $Q_r=\binom dr p^r(1-p)^{d-r}$. Then
\begin{equation}\label{eq:joint-binomial}
\Pp(|A\cap R|=|B\cap R|=r)
\le \exp\left(\frac{2b}{d}+\frac{4br^2}{d}\right)Q_r^2.
\end{equation}
\end{lemma}

\begin{proof}
Summing over the possible values of $|A\cap B\cap R|$ and using the independence of the selections in $A\cap B$, $A\setminus B$, and $B\setminus A$, we obtain
\[
\begin{aligned}
\Pp(|A\cap R|=|B\cap R|=r)
&=\sum_{i=0}^{\min\{b,r\}}\binom bi p^i(1-p)^{b-i}
\left(\binom{d-b}{r-i}p^{r-i}(1-p)^{d-b-r+i}\right)^2\\
&=\sum_{i=0}^{\min\{b,r\}}\binom bi\binom{d-b}{r-i}^{\!2}p^{2r-i}(1-p)^{2d-b-2r+i}.
\end{aligned}
\]
Dividing by $Q_r^2$ gives
\begin{equation}\label{eq:joint-ratio}
\sum_{i=0}^{\min\{b,r\}}\binom bi
\left(\frac{\binom{d-b}{r-i}}{\binom dr}\right)^{\!2}
p^{-i}(1-p)^{-b+i}.
\end{equation}
For every $0\le i\le\min\{b,r\}$, we have
$
\frac{\binom{d-b}{r-i}}{\binom dr}
\le\frac{\binom{d-i}{r-i}}{\binom dr}
=\prod_{j=0}^{i-1}\frac{r-j}{d-j}
\le\left(\frac{2r}{d}\right)^i,
$
where the last inequality follows from $r\le d/2$. Since $\binom bi\le b^i/i!$ and $p=1/d$, the expression in~\eqref{eq:joint-ratio} is at most
$
(1-p)^{-b}\sum_{i=0}^{\infty}\frac{(4br^2/d)^i}{i!}
=(1-p)^{-b}\exp\left(\frac{4br^2}{d}\right).
$
The desired bound follows from $(1-1/d)^{-b}\le\exp(2b/d)$.
\end{proof}

We also use the standard consequence of Tur\'an's theorem that $\alpha(H)\ge |V(H)|/(\overline d(H)+1)$ for every nonempty graph $H$; see, for example, \cite[proof of Theorem~2]{Glock}.

%\subsection{Proof of the edge-codegree theorem}

%We prove the more explicit conclusion $\overline d(H)\ge\eta_0 k(d)$, where
%\begin{equation}\label{eq:edge-constant}
% \eta_0=\frac{1}{25(160000e^5+200e^2+1)}.
%\end{equation}
%Theorem~\ref{thm:edge-codegree} then holds with $c_*=\eta_0/2$.

%\subsection{Proof of the edge-codegree theorem}

%We prove that $G$ contains an induced bipartite subgraph $H$ with $\overline d(H)\ge\eta_0 k(d)$, where $\eta_0=1/(25(160000e^5+200e^2+1))$. Since $k(d)\ge\frac12\log d/\log\log d$ for $d\ge3$, this gives Theorem~\ref{thm:edge-codegree} with $c_*=\eta_0/2$.

\begin{proof}[Proof of Theorem~\ref{thm:edge-codegree}]
Let $k=k(d)$ and $\eta_0=1/[25(160000e^5+200e^2+1)]$, and choose $c=\eta_0/2$. For $3\le d\le15$, we have $c\log d/\log\log d<1$, so any edge of $G$ gives the desired induced bipartite subgraph. We may therefore assume $d\ge16$, in which case $2\le k\le d/2$.

By Lemma~\ref{lem:degenerate-reduction}, we may assume that $G$ is $d$-degenerate. Fix a $d$-degeneracy ordering, and for each vertex $v$ choose $N_v$ with $\NL(v)\subseteq N_v\subseteq N_G(v)$ and $|N_v|=d$. Write $n=|V(G)|$ and $p=1/d$.

Let $X$ be a random subset of $V(G)$ obtained by selecting each vertex independently with probability $p$. Define $S=\{x\in X:\NL(x)\cap X=\varnothing\}$ and $Y=\{y\notin X:|N_y\cap X|=k\}$. The set $S$ is independent. We also consider $Y_0=\{y\in Y:|N_y\cap S|\ge k/50\}$. 
Since $|N_y|=d$ and $y\notin N_y$, we have $\Pp(y\in Y)=q:=(1-p)\binom dk p^k(1-p)^{d-k}$ for every $y\in V(G)$. Thus, $\E|X|=pn$ and $\E|Y|=qn$. Since $k\le\log d/\log\log d$ and $k\le\log d$, we have $k^k\le d$. Using $\binom dk\ge(d/k)^k$ and $(1-1/d)^d\ge e^{-2}$, we obtain
\begin{equation}\label{eq:edge-q-lower}
q\ge k^{-k}(1-p)^d\ge e^{-2}p.
\end{equation}

For every edge $uv$, we have $|N_u\cap N_v|\le |N_G(u)\cap N_G(v)|\le d/k^2$. By Lemma~\ref{lem:joint-binomial},
\[
\begin{aligned}
\Pp(u,v\in Y)
\le\Pp(|N_u\cap X|=|N_v\cap X|=k)
\le\exp\left(4+\frac{2}{k^2}\right)
\left(\frac{q}{1-p}\right)^2
\le2e^5q^2.
\end{aligned}
\]
Here we used $q=(1-p)\Pp(\operatorname{Bin}(d,p)=k)$, $2/k^2\le1$, and $(1-p)^{-2}\le(16/15)^2<2$.

We next show that $\E|Y_0|\ge qn/50$. Fix $y\in V(G)$ and condition on $y\in Y$. Then $K=N_y\cap X$ is a uniformly random $k$-subset of $N_y$. Fix $x\in N_y$ and write $h_x=|\NL(x)\cap N_y|$. Since $xy\in E(G)$, $h_x\le d/k^2$. Conditional on $y\in Y$ and $x\in K$, the set $K\setminus\{x\}$ is uniformly distributed over the $(k-1)$-subsets of $N_y\setminus\{x\}$. Hence
\[
\begin{aligned}
\Pp(K\cap\NL(x)=\varnothing\mid y\in Y,\,x\in K)
&=\frac{\binom{d-1-h_x}{k-1}}{\binom{d-1}{k-1}}\\
&\ge\left(1-\frac{2h_x}{d}\right)^{k-1}
\ge\exp\left(-\frac{4h_x(k-1)}{d}\right)\ge e^{-1}.
\end{aligned}
\]
Here we used $k\le d/2$ and $h_x\le d/k^2$. The random choices in $\NL(x)\setminus(N_y\cup\{y\})$ remain independent of the conditioning, and there are at most $d$ such vertices. Since $y\notin X$, we obtain $\Pp(x\in S\mid y\in Y,\,x\in K)\ge e^{-1}(1-p)^d\ge e^{-3}>1/25$.
It follows that $\E(|K\setminus S|\mid y\in Y)\le24k/25$. By Markov's inequality,
$
\Pp(y\in Y_0\mid y\in Y)\ge1-\frac{24k/25}{49k/50}=\frac1{49}>\frac1{50}.
$
Summing over all vertices yields $\E|Y_0|\ge\frac{qn}{50}.$

By the preceding estimates, $\E|Y_0|\ge qn/50$, $\E e(G[Y])\le2e^5ndq^2$, and $\E|X|=pn$. Hence
\[
\begin{aligned}
&\E\left[|Y_0|-\frac{e(G[Y])}{400e^5qd}
-\frac{q|X|}{200p}\right]
\ge\frac{qn}{50}
-\frac{2e^5ndq^2}{400e^5qd}
-\frac{qpn}{200p}
=\frac{qn}{100}>0.
\end{aligned}
\]
Thus, there exists a choice of $X$ such that
$
|Y_0|>\frac{e(G[Y])}{400e^5qd}+\frac{q|X|}{200p}.
$
Then $Y_0\ne\varnothing$, $e(G[Y_0])\le e(G[Y])<400e^5qd|Y_0|$, and $|S|\le|X|<200p|Y_0|/q$. In particular, $G[Y_0]$ has average degree less than $800e^5qd$. By Tur\'an's theorem, $Y_0$ contains an independent set $T$ with $|T|\ge|Y_0|/(800e^5qd+1)$. Hence, using $pd=1$ and~\eqref{eq:edge-q-lower},
$
|S|\le|X|<\frac{200p}{q}(800e^5qd+1)|T|\le(160000e^5+200e^2)|T|.
$
Since $S\subseteq X$ and $T\subseteq Y_0\subseteq V(G)\setminus X$, the sets $S$ and $T$ are disjoint and independent. Every vertex of $T$ has at least $k/50$ neighbors in $S$, so $H=G[S\cup T]$ is an induced bipartite subgraph satisfying
\[
\overline d(H)\ge\frac{2(k/50)|T|}{|S|+|T|}\ge\frac{k}{25(160000e^5+200e^2+1)}=\eta_0 k.
\]
Since $k(d)\ge\frac12\log d/\log\log d$, we obtain $\overline d(H)\ge c\log d/\log\log d$, as required.
\end{proof}

\section{Concluding remarks}\label{sec:questions}

In this paper, we studied dense induced bipartite subgraphs under local sparsity. The methods developed in this paper also allow us to extend Theorem~\ref{thm:main} to $K_{s,t}$-free graphs and Theorem~\ref{thm:edge-codegree} to a more general edge codegree condition. The corresponding results are stated below.

\begin{theorem}\label{thm:edge-codegree-parameter}
There exist an absolute constant $c>0$ and an absolute integer $d_0$ such that the following holds. Let $d\ge d_0$ and $2\le r\le\lfloor\log d/\log\log d\rfloor$ be integers. If $G$ has minimum degree at least $d$ and every edge belongs to at most $d/r^2$ triangles, then $G$ contains an induced bipartite subgraph of average degree at least $cr$.
\end{theorem}

The proof of Theorem~\ref{thm:edge-codegree} applies with $r$ in place of $k(d)$. Indeed, $r^r\le d$ gives the lower bound $q\ge e^{-2}/d$ in~\eqref{eq:edge-q-lower}, while the edge-codegree assumption bounds the exponential factor in Lemma~\ref{lem:joint-binomial} by $\exp(2/r^2+4)$. The estimate for the number of sampled neighbors that remain in $S$ and all subsequent bounds hold with the same absolute constants.

\begin{theorem}\label{thm:asymmetric-fractional}
There exists an absolute integer $d_0$ such that the following holds. Let $d\ge d_0$ and $2\le s\le t$ be integers satisfying $\log d-\log((t-1)/(s-1))-16s\log\log d>0$. Every $d$-degenerate $K_{s,t}$-free graph $G$ satisfies
\[
\chif(G)\le\frac{d}{\log d-\log\!\left(\frac{t-1}{s-1}\right)-16s\log\log d}.
\]
\end{theorem}

To obtain this extension, fix a target vertex $v$ and let $A=N_L(v)$. Any $s-1$ vertices other than $v$ have at most $t-1$ common neighbors in $A$. For $s\ge3$, the argument of Lemma~\ref{lem:weighted}, with its parameter $s$ replaced by $s-1$, therefore bounds the number of vertices with weighted neighborhood in $A$ greater than $\delta$ by $2(s-2)M/\delta$, provided that $\delta\ge2(t-1)^{1/(s-1)}h^{1/(s-1)}M^{1-1/(s-1)}$. Write $L=\log d$ and $\lambda=\log((t-1)/(s-1))$, and choose $h=\frac{s-1}{t-1}L^{-16(s-1)}$, $R=L^4$, $\delta=L^{-8}$, and $\alpha=(L-\lambda-(16s-13)\log L)/d$. The estimates in the proof of Theorem~\ref{thm:main} then apply uniformly under the stated positivity condition. For $s=2$, every vertex other than $v$ has at most $t-1$ neighbors in $A$, so $Y_i\le(t-1)h=L^{-16}<\delta$, and all steps are light. In particular, this extension allows fixed $s$ and $t=d^{o(1)}$ while retaining $\chif(G)\le(1+o(1))d/\log d$.

Returning to the edge-codegree condition, a natural question is whether the bound in
Theorem~\ref{thm:edge-codegree} can be relaxed by a factor of order
$\log d/\log\log d$.

\begin{question}\label{question:edge}
Do there exist absolute constants $c>0$ and $d_0$ such that every
graph $G$ with minimum degree at least $d\ge d_0$ and
$c_G(uv)\le d\log\log d/(2\log d)$ for every edge $uv$ contains
an induced bipartite subgraph of average degree at least
$c\log d/\log\log d$?
\end{question}

We now restate a conjecture of Martinsson and Steiner concerning the fractional chromatic number of triangle-free graphs. Our results establish the asymptotically optimal leading constant for $K_{t,t}$-free graphs. Moreover, Theorem~\ref{thm:triangle-free} improves the leading constant $4$ in the general bound of Martinsson and Steiner to a constant strictly smaller than $3$. We conjecture the optimal leading constant $1$ can be achieved for all triangle-free graphs.

\begin{conjecture}
Every triangle-free $d$-degenerate graph $G$ satisfies
\[
\chi_f(G)\le (1+o(1))\frac{d}{\log d}.
\]
\end{conjecture}

Allen, Dhawan, and Noel~\cite[Problem~4.1]{ADN} asked whether every $d$-degenerate graph admitting a degeneracy ordering with $|\NL(u)\cap\NR(v)|\le d^{1-o(1)}$ for all vertices $u,v$ satisfies $\chif(G)\le(1+o(1))d/\log d$. The following construction gives a negative answer to this question as stated.

\begin{proposition}\label{prop:ordered-codegree-counterexample}
For infinitely many integers $d$, there exists a $d$-degenerate graph $G$ admitting a $d$-degeneracy ordering such that $|\NL(x)\cap\NR(y)|\le2\sqrt d$ for all $x,y\in V(G)$, but $\chif(G)\ge(2-o(1))d/\log d$.
\end{proposition}

\begin{proof}
By~\cite[Theorem~1.4]{ADN}, for every sufficiently large integer $q$ there exists a $q$-degenerate graph $H=H_q$ of girth at least $5$ with $\chif(H)\ge(1-o(1))q/\log q$. Replace each vertex $v$ of $H$ by a clique $B_v$ of order $q$, and join every vertex of $B_u$ to every vertex of $B_v$ whenever $uv\in E(H)$. Let $G$ be the resulting graph and let $d=q^2+q-1$.

Fix a $q$-degeneracy ordering of $H$ and order the vertices of $G$ by placing the blocks $B_v$ consecutively in this order, with an arbitrary order within each block. Every vertex has at most $q^2$ left-neighbors in earlier blocks and at most $q-1$ in its own block. Hence this is a $d$-degeneracy ordering of $G$.

We claim that $|\NL(x)\cap\NR(y)|\le2q-2$ for all $x,y\in V(G)$. Any vertex in this intersection must lie strictly between $y$ and $x$ in the ordering. If $x,y$ belong to the same block, the intersection is contained in that block and has size at most $q-2$. Suppose that $x\in B_u$ and $y\in B_v$ with $u\ne v$. If $uv\in E(H)$, then triangle-freeness of $H$ implies that every common neighbor of $x$ and $y$ belongs to $B_u\cup B_v$, giving at most $2q-2$ common neighbors. If $uv\notin E(H)$, then $C_4$-freeness of $H$ implies that $u$ and $v$ have at most one common neighbor, so $x$ and $y$ have at most $q$ common neighbors. This proves the claim. Since $2q-2\le2\sqrt d$, the required intersection bound follows.

To bound $\chif(G)$ from below, choose nonnegative vertex weights $w_v$ on $H$ such that $\sum_{v\in V(H)}w_v=\chif(H)$ and $\sum_{v\in J}w_v\le1$ for every independent set $J$ of $H$. Assign weight $w_v$ to each vertex of $B_v$. Every independent set of $G$ meets each block in at most one vertex, and the blocks it meets correspond to an independent set of $H$. Thus these weights form a feasible fractional clique of $G$, and consequently
\[
\chif(G)\ge q\chif(H)\ge(1-o(1))\frac{q^2}{\log q}
=(2-o(1))\frac{d}{\log d}.
\]
\end{proof}

For every function $\eta(d)\to0$, we have $2\sqrt d\le d^{1-\eta(d)}$ for all sufficiently large $d$. Thus Proposition~\ref{prop:ordered-codegree-counterexample} excludes the proposed conclusion even under a substantially stronger intersection bound. The construction above contains triangles. We conjecture that the upper bound proposed in~\cite[Problem~4.1]{ADN} holds under the additional assumption that $G$ is triangle-free.

\begin{conjecture}\label{conj:triangle-free-ordered-codegree}
Let $G$ be a triangle-free $d$-degenerate graph admitting a $d$-degeneracy ordering such that $|\NL(u)\cap\NR(v)|\le d^{1-o(1)}$ for all $u,v\in V(G)$. Then $\chif(G)\le(1+o(1))d/\log d$.
\end{conjecture}

\section*{Declaration of AI use}
During the preparation of this manuscript, the authors used ChatGPT (OpenAI, including GPT-5.4 and GPT-6 Pro) to assist with mathematical discussions, proof development and verification, and language editing.

The research problems, main objectives, and overall mathematical framework were developed by the authors. The results on \(K_{t,t}\)-free graphs and bounded edge-codegrees, including their applications to dense induced bipartite subgraphs, were developed primarily by the authors, with AI assistance in reviewing arguments and improving the exposition.

ChatGPT contributed to the development and refinement of the iterative entropy estimate in Lemma \ref{lem:tf-record-entropy}, the analysis of the stopped weight processes, and the optimization of constants leading to the bound \((2.85+o(1))d/\log d\). The authors selected, examined, and revised these arguments and integrated them into the manuscript.

The authors take full responsibility for the correctness of all mathematical results and proofs, as well as the final content of the manuscript.

\end{document}